\documentclass[11pt]{article}
\usepackage{bbm}
\usepackage{mathrsfs}
\usepackage{amscd}
\usepackage{amsmath,amsfonts,amssymb,amscd}
\usepackage{indentfirst,graphics,epsfig,psfrag}
\input{epsf}
\usepackage{ifpdf}
\usepackage{enumerate}
\usepackage{appendix}
\usepackage{enumerate}
\usepackage{color}
\usepackage{lineno}
\usepackage{authblk}
\makeatletter \@addtoreset{figure}{section} \makeatother
\makeatletter
\long\def\@makecaption#1#2{%
   \vskip 10\p@
   \setbox\@tempboxa\hbox{{#1}\ \ #2}%
   \ifdim \wd\@tempboxa >\hsize

       {#1}\ \ #2\par
   \else
       \hbox to\hsize{\hfil\box\@tempboxa\hfil}%
   \fi}
\makeatother

\newtheorem{thm}{Theorem}
\newtheorem{cor}{Corollary}
\newtheorem{lem}{Lemma}

\newtheorem{obs}{Observation}
\newtheorem{pro}{Proposition}

\newcommand{\qed}{{\hfill\rule{3pt}{7pt}}}

\def\qed{\hfill \rule{4pt}{7pt}}

\begin{document}
\title{\textbf{The Steiner $(n-3)$-diameter of a graph}
\footnote{Supported by the National Science Foundation of China
(Nos. 11161037, 11101232, and 11461054) and the Science Found of Qinghai Province
(Nos. 2016-ZJ-948Q, and 2014-ZJ-907).}}
\author[1]{Yaping Mao\footnote{E-mail: maoyaping@ymail.com}}
\author[2]{Christopher Melekian}
\author[2]{Eddie Cheng}
\affil[1]{Department of Mathematics, Qinghai Normal University}
\affil[2]{Department of Mathematics and Statistics, Oakland University}

\date{}
\maketitle
\begin{abstract}
The Steiner distance of a graph, introduced by Chartrand,
Oellermann, Tian and Zou in 1989, is a natural generalization of the
concept of classical graph distance. For a connected graph $G$ of
order at least $2$ and $S\subseteq V(G)$, the \emph{Steiner
distance} $d(S)$ among the vertices of $S$ is the minimum size among
all connected subgraphs whose vertex sets contain $S$. Let $n$ and
$k$ be two integers with $2\leq k\leq n$. Then the \emph{Steiner
$k$-eccentricity $e_k(v)$} of a vertex $v$ of $G$ is defined by
$e_k(v)=\max \{d(S)\,|\,S\subseteq V(G), \ |S|=k, \ and \ v\in S
\}$. Furthermore, the Steiner \emph{$k$-diameter} of $G$ is
$sdiam_k(G)=\max \{e_k(v)\,|\, v\in V(G)\}$. In 2011, Chartrand,
Okamoto, Zhang showed that $k-1\leq sdiam_k(G)\leq n-1$. In this
paper, graphs with $sdiam_k(G)=\ell$ for $k=n,n-1,n-2,n-3$ and
$k-1\leq \ell \leq n-1$
are characterized, respectively. \\[2mm]
\end{abstract}

\section{Introduction}

All graphs in this paper are undirected, finite and simple. We refer
to \cite{Bondy} for graph-theoretic notation and terminology not
described here. We
divide our introduction into the following four subsections to state
the motivations and our results of this paper.

\subsection{Distance and its generalizations}

Distance is one of the most basic concepts of
graph theory. If $G$ is a connected graph and $u,v\in
V(G)$, then the \emph{distance} $d(u,v)$ between $u$ and $v$ is the
length of a shortest path connecting $u$ and $v$. If $v$ is a vertex
of a connected graph $G$, then the \emph{eccentricity} $e(v)$ of $v$
is defined by $e(v)=\max\{d(u,v)\,|\,u\in V(G)\}$. Furthermore, the
\emph{radius} $rad(G)$ and \emph{diameter} $diam(G)$ of $G$ are
defined by $rad(G)=\min\{e(v)\,|\,v\in V(G)\}$ and $diam(G)=\max
\{e(v)\,|\,v\in V(G)\}$. These last two concepts are related by the
inequalities $rad(G)\leq diam(G) \leq 2 rad(G)$. The \emph{center}
$C(G)$ of a connected graph $G$ is the subgraph induced by the
vertices $u$ of $G$ with $e(u)=rad(G)$. Recently, Goddard and
Oellermann gave a survey paper on this subject, see \cite{Goddard}.

The distance between two vertices $u$ and $v$ in a connected graph
$G$ also equals the minimum size of a connected subgraph of $G$
containing both $u$ and $v$. This observation suggests a
generalization of distance. The Steiner distance of a graph,
introduced by Chartrand, Oellermann, Tian and Zou in 1989, is a
natural and nice generalization of the concept of classical graph
distance. For a graph $G(V,E)$ and a set $S\subseteq V(G)$ of at
least two vertices, \emph{an $S$-Steiner tree} or \emph{a Steiner
tree connecting $S$} (or simply, \emph{an $S$-tree}) is a such
subgraph $T(V',E')$ of $G$ that is a tree with $S\subseteq V'$. Let
$G$ be a connected graph of order at least $2$ and let $S$ be a
nonempty set of vertices of $G$. Then the \emph{Steiner distance}
$d_G(S)$ among the vertices of $S$ (or simply the distance of $S$)
is the minimum size among all connected subgraphs whose vertex sets
contain $S$. Note that if $H$ is a connected subgraph of $G$ such
that $S\subseteq V(H)$ and $|E(H)|=d_G(S)$, then $H$ is a tree.
Observe that $d_G(S)=\min\{e(T)\,|\,S\subseteq V(T)\}$, where $T$ is
subtree of $G$. Furthermore, if $S=\{u,v\}$, then $d_G(S)=d(u,v)$ is
nothing new but the classical distance between $u$ and $v$. Set
$d_G(S)=\infty$ when there is no $S$-Steiner tree in $G$.

\begin{obs}\label{obs1}
Let $G$ be a graph of order $n$ and $k$ be an integer with $2\leq k\leq n$. If $S\subseteq V(G)$ and $|S|=k$, then $d_G(S)\geq k-1$.
\end{obs}



Let $n$ and $k$ be two integers with $2\leq k\leq n$. The
\emph{Steiner $k$-eccentricity $e_k(v)$} of a vertex $v$ of $G$ is
defined by $e_k(v)=\max \{d(S)\,|\,S\subseteq V(G), |S|=k,~and~v\in
S \}$. The \emph{Steiner $k$-radius} of $G$ is $srad_k(G)=\min \{
e_k(v)\,|\,v\in V(G)\}$, while the \emph{Steiner $k$-diameter} of
$G$ is $sdiam_k(G)=\max \{e_k(v)\,|\,v\in V(G)\}$. Note for every
connected graph $G$ that $e_2(v)=e(v)$ for all vertices $v$ of $G$
and that $srad_2(G)=rad(G)$ and $sdiam_2(G)=diam(G)$.


\begin{obs}\label{obs2}
Let $k,n$ be two integers with $2\leq k\leq n$.

$(1)$ If $H$ is a spanning subgraph of $G$, then $sdiam_k(G)\leq
sdiam_k(H)$.

$(2)$ For a connected graph $G$, $sdiam_k(G)\leq sdiam_{k+1}(G)$.
\end{obs}

As a generalization of the center of a graph, the \emph{Steiner
$k$-center} $C_k(G)\ (k\geq 2)$ of a connected graph $G$ is the
subgraph induced by the vertices $v$ of $G$ with $e_k(v)=srad_k(G)$.
Oellermann and Tian \cite{OellermannT} showed that every graph is
the $k$-center of some graph. In particular, they showed that the
$k$-center of a tree is a tree and those trees that are $k$-centers
of trees are characterized. The \emph{Steiner $k$-median} of $G$ is
the subgraph of $G$ induced by the vertices of $G$ of minimum
Steiner $k$-distance. For Steiner centers and Steiner medians, we
refer to \cite{Oellermann, Oellermann2, OellermannT}.

The \emph{average Steiner distance} $\mu_k(G)$ of a graph $G$,
introduced by Dankelmann, Oellermann and Swart in
\cite{DankelmannOS}, is defined as the average of the Steiner
distances of all $k$-subsets of $V(G)$, i.e.
$$
\mu_k(G)={n\choose k}^{-1}\sum_{S\subseteq V(G),|S|=k}d_G(S).
$$
For more details on average Steiner distance, we refer to
\cite{DankelmannOS, DankelmannSO}.

Let $G$ be a $k$-connected graph and $u$, $v$ be any pair of
vertices of $G$. Let $P_k(u,v)$ be a family of $k$ vertex-disjoint
paths between $u$ and $v$, i.e., $P_k(u,v)=\{P_1,P_2,\cdots,P_k\}$,
where $p_1\leq p_2\leq \cdots \leq p_k$ and $p_i$ denotes the number
of edges of path $P_i$. The \emph{$k$-distance} $d_k(u,v)$ between
vertices $u$ and $v$ is the minimum $|p_k|$ among all $P_k(u,v)$ and
the \emph{$k$-diameter} $d_k(G)$ of $G$ is defined as the maximum
$k$-distance $d_k(u,v)$ over all pairs $u,v$ of vertices of $G$. The
concept of $k$-diameter emerges rather naturally when one looks at
the performance of routing algorithms. Its applications to network
routing in distributed and parallel processing are studied and
discussed by various authors including Chung \cite{Chung}, Du, Lyuu
and Hsu \cite{Du}, Hsu \cite{Hsu, Hsu2}, Meyer and Pradhan
\cite{Meyer}.

\subsection{Application background and progress of Steiner distance}

The Steiner tree problem in networks, and particularly in graphs,
was formulated quite recently-in 1971-by Hakimi (see \cite{Hakimi})
and Levi (see \cite{Levi}). In the case of an unweighted, undirected
graph, this problem consists of finding, for a subset of vertices
$S$, a minimal-size connected subgraph that contains the vertices in
$S$. The computational side of this problem has been widely studied,
and it is known that it is an NP-hard problem for general graphs
(see \cite{HwangRW}). The determination of a Steiner tree in a graph
is a discrete analogue of the well-known geometric Steiner problem:
In a Euclidean space (usually a Euclidean plane) find the shortest
possible network of line segments interconnecting a set of given
points. Steiner trees have application to multiprocessor computer
networks. For example, it may be desired to connect a certain set of
processors with a subnetwork that uses the least number of
communication links. A Steiner tree for the vertices, corresponding
to the processors that need to be connected, corresponds to such a
desired subnetwork.

In \cite{ChartrandOZ}, Chartrand, Okamoto, Zhang obtained the
following result.

\begin{thm}{\upshape\cite{ChartrandOZ}}\label{th1}
Let $k,n$ be two integers with $2\leq k\leq n$, and let $G$ be a
connected graph of order $n$. Then $k-1\leq sdiam_k(G)\leq n-1$.
Moreover, the upper and lower bounds are sharp.
\end{thm}

In \cite{DankelmannSO2}, Dankelmann, Swart and Oellermann obtained a
bound on $sdiam_k(G)$ for a graph $G$ in terms of the order of $G$
and the minimum degree of $G$, that is, $sdiam_k(G)\leq
\frac{3p}{\delta+1}+3n$. Later, Ali, Dankelmann, Mukwembi
\cite{AliDM} improved the bound of $sdiam_k(G)$ and showed that
$sdiam_k(G)\leq \frac{3p}{\delta+1}+2n-5$ for all connected graphs
$G$. Moreover, they constructed graphs to show that the bounds are
asymptotically best possible.

In \cite{Bloom}, Bloom characterized graphs with diameter $2$. Mao
\cite{Mao} characterized the graphs with $sdiam_3(G)=2,3,n-1$, and
obtained the Nordhaus-Gaddum-type results for the parameter
$sdiam_k(G)$. In this paper, graphs with $sdiam_k(G)=\ell$ for
$k=n,n-1,n-2,n-3$ and $k-1\leq \ell \leq n-1$ are characterized,
respectively.

\section{The case $k=n,n-1,n-2$}

To begin with, we show the following two lemmas, which will be used
later.

\begin{lem}\label{lem1}
Let $k,n$ be two integers with $2\leq k\leq n$, and let $T$ be a
tree of order $n$. Then $sdiam_k(T)=n-1$ if and only if $r\leq k$,
where $r$ is the number of the leaves in $T$.
\end{lem}
\begin{pf}
Suppose $r\leq k$. Let $v_1,v_2,\cdots,v_{r}$ be all the leaves of
$T$. Choose $S\subseteq V(T)$ and $|S|=k$ such that
$v_1,v_2,\cdots,v_{r}\in S$. Then any tree connecting $S$ must use
all edges of $T$. Since $|E(T)|=n-1$, it follows that $d_{T}(S)\geq
|E(T)|=n-1$. From the arbitrariness of $S$, we have $sdiam_k(T)\geq
n-1$. Combining this with Theorem \ref{th1}, we have
$sdiam_k(T)=n-1$.

Conversely, suppose $sdiam_k(T)=n-1$. If $s\geq k+1$, then for any
$S\subseteq V(G)$ with $|S|=k$, there exists a leaf $v\in V(T)$ such
that $v\notin S$. Set $T'=T\setminus v$. Then the tree $T'$ is an
$S$-Steiner tree and hence $d_{T}(S)\leq n-2$. From the
arbitrariness of $S$, we have $sdiam_k(T)\leq n-2<n-1$, a
contradiction. Therefore, $s\leq k$.\qed
\end{pf}

\begin{lem}\label{lem2}
Let $k,n$ be two integers with $1\leq k\leq n-2$. Let $G$ be a
connected graph of order $n$. Then $sdiam_{n-k}(G)=n-1$ if and only if $G$ contains at least $k$ cut vertices.

\end{lem}
\begin{pf}
Suppose $G$ contains at least $k$ cut vertices; pick $k$ of them, say $v_1,v_2,\cdots,v_k$. Choose
$S=V(G)\setminus \{v_1,v_2,\cdots,v_k\}$. Then $|S|=n-k$ and any
$S$-Steiner tree $T$ must occupy the vertices $v_1,v_2,\cdots,v_k$,
which implies that $|V(T)|=n$ and $e(T)=n-1$. Furthermore,
$d_{G}(S)\geq e(T)=n-1$ and hence $sdiam_{n-k}(G)\geq d_{G}(S)\geq
n-1$. Theorem \ref{th1} yields $sdiam_{n-k}(G)=n-1$, as desired.

Conversely, suppose that $sdiam_{n-k}(G) = n-1$.
Assume to the contrary that $G$ contains at most $k-1$ cut vertices; let $C$ be the set of all cut vertices in $G$.
Then for any $S \subset V(G)$ with $|S| = n-k$, we have $|S \cup C| \leq n-1$,
so we can find a vertex $x \in V(G)$ such that $x$ is not a member of $S$ and not a cut vertex of $G$.
Therefore $G \setminus x$ is connected and has a spanning tree $T$.
Observe that $|V(T)| = n-1$, so $d_G(S)\leq |e(T)| = n-2$.
Since $S$ was arbitrary, we have $sdiam_{n-k}(G) \leq n-2$, a contradiction. So $G$ contains at least $k$ cut vertices, as desired.
\qed
\end{pf}

\begin{pro}\label{pro1}
Let $k,n$ be two integers with $1\leq k\leq n-2$, and let $G$ be a
graph of order $n$. Then $\kappa(G)\geq k$ if and only if
$sdiam_{n-k+1}(G)=n-k$.
\end{pro}
\begin{pf}
For any $S\subseteq V(G)$ with $|S|=n-k+1$, we have $|V(G)\setminus
S|=k-1$. Since $\kappa(G)\geq k$, it follows that $G[S]$ is
connected. Therefore, $G[S]$ contains a spanning tree $T$ of order
$n-k+1$ and hence $e(T)=n-k$. From the arbitrariness of $S$, we have
$sdiam_{n-k+1}(G)\leq d_{T}(S)=e(T)=n-k$. From this together with
Theorem \ref{th1}, $sdiam_{n-k+1}(G)=n-k$.

Conversely, we suppose $sdiam_{n-k+1}(G)=n-k$. If $\kappa(G)\leq
k-1$, then there exist a cut set $U\subseteq V(G)$ with
$|U|=\kappa(G)$ such that $G\setminus U$ is disconnected. Let
$C_1,C_2,\cdots,C_r$ be the connected components of $G\setminus U$.
Note that $(\bigcup_{i=1}^{r}V(C_i))\bigcup U=V(G)$. Since $|U|\leq
k-1$, it follows that $|\bigcup_{i=1}^{r}V(C_i)|\geq n-(k-1)=n-k+1$.
Pick up $n-k+1$ vertices from $\bigcup_{i=1}^{r}V(C_i)$. Let $S$ be
the vertex set of these $n-k+1$ vertices. Then any $S$-Steiner tree
$T$ must use at least one vertex of $U$, which implies that
$|V(T)|\geq (n-k+1)+1=n-k+2$. Thus, $d_G(S)\geq e(T)=|V(T)|-1\geq
n-k+1$ and hence $sdiam_{n-k+1}(G)\geq d_G(S)\geq n-k+1$, a
contradiction. So $G$ is $k$-connected.\qed \vskip 0.3em
\end{pf}

For $k=n,n-1$, we have $sdiam_n(G)=n-1$ and $n-2\leq
sdiam_{n-1}(G)\leq n-1$ by Theorem \ref{th1}. The following two
corollaries are immediate by Proposition \ref{pro1}.

\begin{cor}\label{cor1}
Let $G$ be a graph of order $n$. Then $sdiam_n(G)=n-1$ if and only
if $G$ is connected.
\end{cor}

\begin{cor}\label{cor2}
Let $G$ be a connected graph of order $n$. Then

$(1)$ $sdiam_{n-1}(G)=n-2$ if and only if $G$ is $2$-connected;

$(2)$ $sdiam_{n-1}(G)=n-1$ if and only if $G$ contains at least one
cut vertex.
\end{cor}

From Theorem \ref{th1}, we know $n-3\leq sdiam_{n-2}(G)\leq n-1$.
Let us now characterize the graphs with
$sdiam_{n-2}(G)=n-3,n-2,n-1$.

\begin{thm}\label{th2}
Let $G$ be a connected graph of order $n \ (n\geq 4)$. Then

$(1)$ $sdiam_{n-2}(G)=n-3$ if and only if $\kappa(G)\geq 3$.

$(2)$ $sdiam_{n-2}(G)=n-2$ if and only if $\kappa(G)=2$ or $G$
contains only one cut vertex.

$(3)$ $sdiam_{n-2}(G)=n-1$ if and only if there are at least two cut
vertices in $G$.
\end{thm}
\begin{pf}
$(1)$ The result follows by Proposition \ref{pro1}.

$(3)$ The result follows by Lemma \ref{lem2}.

$(2)$ By $(1)$, we must have $\kappa(G) \leq 2$. If $\kappa(G) = 2$, we are done. If $\kappa(G) = 1$, then $G$ must contain at least one cut vertex, but fewer than two cut vertices by $(3)$, so $G$ contains exactly one cut vertex. \qed

\end{pf}

\section{The case $k=n-3$}

From Theorem \ref{th1}, we know that $n-4\leq sdiam_{n-3}(G)\leq
n-1$. Graphs with $sdiam_{n-3}(G)=n-4,n-3,n-2,n-1$ are characterized
in this section.

The following is immediate from Proposition \ref{pro1} and Lemma \ref{lem2}.

\begin{pro}\label{pro2}
Let $G$ be a connected graph of order $n$. Then $sdiam_{n-3}(G)=n-4$
if and only if $\kappa(G)\geq 4$, and $sdiam_{n-3}(G)=n-1$ if and only if $G$ contains at least $3$ cut vertices.
\end{pro}

\begin{lem}\label{lem3}
Let $G$ be a connected graph of order $n$. If $\kappa(G)=1$, then $sdiam_{n-3}(G)=n-3$
if and only if $G$ satisfies the following two conditions.

$(1)$ $G$ contains only one cut vertex $u$;

$(2)$ for each
connected component $C_i$ of order at least $3$ in $G\setminus u$,
$G[V(C_i)\cup \{u\}]$ is $3$-connected, or $\kappa(G[V(C_i)\cup
\{u\}])=2$ and there exists a vertex $v\in V(C_i)$ such that
$\{u,v\}$ is a vertex cut set of $G[V(C_i)\cup \{u\}]$, and for
each component $C_i^j \ (1\leq j\leq p)$ of $G[V(C_i)\cup
\{u\}]\setminus \{u,v\}$, one of the following conditions holds:

$\bullet$ $uv\in E(G)$;

$\bullet$ $p\geq 3$;

$\bullet$ $p=2$, and $|E_G[v,V(C_i^1)]|\geq 2$ or
$|E_G[v,V(C_i^2)]|\geq 2$.

and one of the following conditions holds:

$\bullet$ $G[V(C_i^j)\cup \{u\}]$ is $3$-connected;

$\bullet$ $G[V(C_i^j)\cup \{v\}]$ is $3$-connected;

$\bullet$ $\kappa(G[V(C_i^j)\cup \{u\}])=\kappa(G[V(C_i^j)\cup
\{v\}])=2$ and $\{y,z\}$ is not a common vertex cut set of
$G[V(C_i^j)\cup \{u\}]$ and $G[V(C_i^j)\cup \{v\}]$ where $z',z''\in
V(C_i^j)$;

$\bullet$ $\kappa(G[V(C_i^j)\cup \{u\}])=2$ and
$\kappa(G[V(C_i^j)\cup \{v\}])=1$ and if $\{z',z''\}$ is a vertex cut
set of $G[V(C_i^j)\cup \{u\}]$, then neither $z'$ nor $z''$ is a cut
vertex of $G[V(C_i^j)\cup \{v\}]$.
\end{lem}
\begin{pf}
In one direction, we suppose $sdiam_{n-3}(G)=n-3$. Assume, to the contrary, that $G$ contains only one cut vertex $u$ such that there
exists a connected component $C_j$ of order at least $3$ in
$G\setminus u$ satisfying one of the following.

$(1)$ $\kappa(G[V(C_j)\cup \{u\}])=1$;

$(2)$ $\kappa(G[V(C_j)\cup \{u\}])=2$ and $\{v,u\}$
is not a vertex cut set of $G[V(C_j)\cup \{u\}]$ for any $v\in V(C_j)$;

$(3)$ $\kappa(G[V(C_j)\cup \{u\}])=2$, $\{v,u\}$ is a
vertex cut set of $G[V(C_j)\cup \{u\}]$, and there exists a component
$C_j^{i'}$ of $G[V(C_j)\cup \{u\}]\setminus \{u,v\}$ satisfying one
of the following conditions.

$\bullet$ $\kappa(G[V(C_j^{i'}\cup \{u\}])=1$.

$\bullet$ $\kappa(G[V(C_j^{i'}\cup \{u\}])=\kappa(G[V(C_j^{i'}\cup
\{v\}])=2$ and $\{y,z\}$ is a common vertex cut set of
$G[V(C_j^{i'}\cup \{u\}]$ and $V(C_j^{i'}\cup \{v\}]$ for any
$y,z\in V(C_i^j)$.

$\bullet$ $\kappa(G[V(C_j^{i'}\cup \{u\}])=2$ and
$\kappa(G[V(C_j^{i'}\cup \{v\}])=1$ and if $\{y,z\}$ is a common
vertex cut set of $G[V(C_j^{i'}\cup \{u\}]$ where $y,z\in V(C_i^j)$,
then at least one of $\{y,z\}$ is a cut vertex of $V(C_j^{i'}\cup
\{v\}]$.

$\bullet$ $uv\notin E(G)$, $p=2$ and there is only one edge between
$v$ and each connected component of $G[V(C_j)\cup \{u\}]\setminus
\{u,v\}$.

Our aim is to show $sdiam_{n-3}(G)\geq n-2$ and get a contradiction. Let
$H_j=G[V(C_j)\cup \{u\}]$. If $\kappa(H_j)=1$, then $u$ is not a cut vertex of $H_j$ since
$C_j$ is a connected component of $G\setminus u$. Therefore, there
exists a cut vertex of $H_j$, say $x$, such that $x\neq u$. Let
$C_j^1,C_j^2,\cdots,C_j^s \ (s\geq 2)$ be the connected components
of $H_j\setminus x$. Clearly, $u\in \bigcup_{i=1}^sV(C_j^i)$.
Without loss of generality, let $u\in V(C_j^1)$. We claim that there
exists a connected component $C_j^{i_1} \ (i_1\in \{2,3,\cdots,s\})$
such that $|E_G[u,C_j^{i_1}]|=0$. Assume, to the contrary, that
$|E_G[u,C_j^{i}]|\geq 1$ for each $i\in \{2,3,\cdots,s\}$. Then
$H_j\setminus x$ is connected, a contradiction. Thus, $x$ is also a
cut vertex in $G$, which contradicts to the fact that $G$ only
contains one cut vertex $u$.

If $\kappa(H_j)=2$ and for any $v\in V(C_j)$ $\{v,u\}$ is not a
vertex cut set of $H_j$. Then there exist a vertex cut $\{x,y\}$ of
$H_j$. Therefore, $H_j\setminus \{x,y\}$ is disconnected. Let
$C_j^1,C_j^2,\cdots,C_j^s \ (s\geq 2)$ be the connected components
of $H_j\setminus \{x,y\}$. Clearly, $u\in \bigcup_{i=1}^sV(C_j^i)$.
Without loss of generality, let $u\in V(C_j^1)$. Then there exists a
connected component $C_j^{i_1} \ (i_1\in \{2,3,\cdots,s\})$ such
that $|E_G[u,C_j^{i_1}]|=0$. Choose $S=V(G)\setminus \{x,y,u\}$.
Then $|S|=n-3$. Obviously, any $S$-Steiner tree, say $T$,
must contain the vertex $u$ and one of $\{x,y\}$, which implies
$|V(T)|\geq n-1$ and hence $d_G(S)\geq e(T)\geq n-2$ and hence
$sdiam_{n-3}(G)\geq d_G(S)\geq n-2$, a contradiction.

Suppose $\kappa(H_j)=2$ and there exists a vertex $v\in V(C_j)$ such
that $\{u,v\}$ is a vertex cut set of $H_j$ and there exists a
component $C_j^{i'}$ of $H_j\setminus \{u,v\}$ such that
$\kappa(G[V(C_j^{i'}\cup \{u\}])=1$. Then there exists a cut vertex
of $G[V(C_j^{i'}\cup \{u\}]$, say $z$. Since $C_j^{i'}$ is
connected, it follows that $z\neq u$. Choose $\bar{S}=\{u,z,v\}$.
Then any $S$-Steiner tree, say $T$, must contain the
vertex $u$ and one of $\{z,v\}$, which implies $|V(T)|\geq n-1$ and
hence $d_G(S)\geq e(T)\geq n-2$ and hence $sdiam_{n-3}(G)\geq
d_G(S)\geq n-2$, a contradiction.

Suppose $\kappa(H_j)=2$ and there exists a vertex $v\in V(C_j)$ such
that $\{u,v\}$ is a vertex cut set of $H_j$ and there exists a
component $C_j^{i'}$ of $H_j\setminus \{u,v\}$ such that
$\kappa(G[V(C_j^{i'}\cup \{u\}])=\kappa(G[V(C_j^{i'}\cup \{v\}])=2$
and $\{y,z\}$ is a common vertex cut set of $G[V(C_j^{i'}\cup
\{u\}]$ and $V(C_j^{i'}\cup \{v\}]$ for any $y,z\in V(C_i^j)$.
Choose $\bar{S}=\{y,z,u\}$. Then any $S$-Steiner tree,
say $T$, must contain the vertex $u$ and one of $\{y,z\}$, which
implies $|V(T)|\geq n-1$ and hence $d_G(S)\geq e(T)\geq n-2$ and
hence $sdiam_{n-3}(G)\geq d_G(S)\geq n-2$, a contradiction.

Suppose $\kappa(H_j)=2$ and there exists a vertex $v\in V(C_j)$ such
that $\{u,v\}$ is a vertex cut set of $H_j$ and there exists a
component $C_j^{i'}$ of $H_j\setminus \{u,v\}$ such that
$\kappa(G[V(C_j^{i'}\cup \{u\}])=2$ and $\kappa(G[V(C_j^{i'}\cup
\{v\}])=1$ and if $\{y,z\}$ is a common vertex cut set of
$G[V(C_j^{i'}\cup \{u\}]$ where $y,z\in V(C_i^j)$, then at least one
of $\{y,z\}$ is a cut vertex of $V(C_j^{i'}\cup \{v\})$. Without
loss of generality, let $y$ is a cut vertex of $V(C_j^{i'}\cup
\{v\}]$. Choose $\bar{S}=\{y,z,u\}$. Then any $S$-Steiner tree, say $T$, must contain the vertex $u$ and one of
$\{y,z\}$, which implies $|V(T)|\geq n-1$ and hence $d_G(S)\geq
e(T)\geq n-2$ and hence $sdiam_{n-3}(G)\geq d_G(S)\geq n-2$, a
contradiction.

Suppose that $uv\notin E(G)$, $p=2$ and there is only one edge
between $v$ and each connected component of $G[V(C_j)\cup
\{u\}]\setminus \{u,v\}$. Let $C_j^1,C_j^2$ be the connected
components of $G[V(C_j)\cup \{u\}]\setminus \{u,v\}$, and let
$vv_p,vv_q$ be the edges between $v$ and $C_j^1,C_j^2$,
respectively. Choose $\bar{S}=\{u,v_p,v_q\}$. Recall that $uv\notin
E(G)$. For any $S$-Steiner tree $T$, $T$ must occupy the
vertex $u$. Also, in order to reach the vertex $v$, $T$ must contain
the vertex one of $\{v_p,v_q\}$, which implies $|V(T)|\geq n-1$ and
hence $d_G(S)\geq e(T)\geq n-2$ and hence $sdiam_{n-3}(G)\geq
d_G(S)\geq n-2$, a contradiction.

Conversely, we suppose that $G$ satisfies the conditions of this
theorem. From Proposition \ref{pro2}, we know that $sdiam_{n-3}(G)\geq
n-3$. So it suffices to show $sdiam_{n-3}(G)\leq n-3$.

In this case, each connected component of $G\setminus u$ is a
connected subgraph of order at least $3$, or an edge of $G$, or an
isolated vertex. Let $w_1,w_2,\cdots,w_r$ be the isolated vertices,
$e_1,e_2,\cdots,e_t$ be the edges, and $C_1,C_2,\cdots,C_r$ be the
connected components of order at least $3$ in $G\setminus u$. Set
$e_i=u_iv_i \ (1\leq i\leq t)$, $W=\{w_1,w_2,\cdots,w_s\}$,
$U=\{u_1,u_2,\cdots,u_t\}$, $V=\{v_1,v_2,\cdots,v_t\}$ and
$n_i=|V(C_i)|$. Obviously, $uw_i,uu_i,uv_i\in E(G)$ and
$s+2t+\sum_{i=1}^rn_i=n-1$. Since $|S|=n-3$, there exists three
vertices $x,y,z$ such that $x,y,z\notin S$ and $x,y,z\in V(G)=W\cup U\cup V\cup
(\cup_{i=1}^rV(C_i))\cup \{u\}$.

Suppose that at least two of $\{x,y,z\}$ belong to $W\cup U\cup V$.
Without loss of generality, let $x,y\in W\cup U\cup V$. Then
$G\setminus \{x,y\}$ is connected and hence $G\setminus \{x,y\}$
contains a spanning tree $T$. Therefore, $d_G(S)\leq e(T)=n-3$, as
desired.

Suppose that only one of $\{x,y,z\}$ belongs to $W\cup U\cup V$.
Without loss of generality, let $x\in W\cup U\cup V$ and $x=w_1$.
Then $y,z\in (\bigcup_{i=1}^rV(C_i))\cup \{u\}$. Without loss of
generality, let $y\in \bigcup_{i=1}^rV(C_i)$. Thus there exists some
$C_j$ such that $y\in C_j$. Since $G[V(C_j)\cup \{u\}]$ is
$2$-connected, it follows that $G[V(C_j)\cup \{u\}]\setminus y$ is
connected and hence $G[V(C_j)\cup \{u\}]\setminus y$ contains a
spanning tree $T_j$. Furthermore, $G[V(C_i)\cup \{u\}]$ contains a
spanning tree $T_i$ for each $i \ (1\leq i\leq r, i\neq j)$. Then
the tree $T$ induced by the edges in $\{uw_i\,|\,2\leq i\leq s\}\cup
\{uu_{i}\,|\,1\leq i\leq t\}\cup \{uv_{i}\,|\,2\leq i\leq t\}\cup
E(T_1)\cup E(T_2)\cup \cdots \cup E(T_r)$ is our desired $S$-Steiner tree. Therefore, $d_G(S)\leq
e(T)=(s-1)+2t+\sum_{i=1,i\neq j}^rn_i+(n_j-1)=n-3$, as desired.

Suppose that none of $\{x,y,z\}$ belongs to $W\cup U\cup V$. Then
$x,y,z\in (\bigcup_{i=1}^rV(C_i))\cup \{u\}$. Therefore, at least
two of $\{x,y,z\}$ belongs to $\bigcup_{i=1}^rV(C_i)$. Without loss
of generality, let $x,y\in \bigcup_{i=1}^rV(C_i)$.

First, we consider the case that $x,y$ belong to different connected
components. Without loss of generality, let $x\in V(C_1)$ and $y\in
V(C_2)$. Since $G[V(C_i)\cup \{u\}] \ (i=1,2)$ is $2$-connected, it
follows that both $G[V(C_1)\cup \{u\}]\setminus x$ and $G[V(C_2)\cup
\{u\}]\setminus y$ is connected. Therefore, $G[V(C_1)\cup
\{u\}]\setminus x$ contains a spanning tree $T_1$ and $G[V(C_2)\cup
\{u\}]\setminus y$ contains a spanning tree $T_2$. Furthermore,
$G[V(C_i)\cup \{u\}]$ contains a spanning tree $T_i$ for each $i \
(3\leq i\leq r)$. Then the tree $T$ induced by the edges in
$\{uw_i\,|\,1\leq i\leq s\}\cup \{uu_{i}\,|\,1\leq i\leq t\}\cup
\{uv_{i}\,|\,1\leq i\leq t\}\cup E(T_1)\cup E(T_2)\cup \cdots \cup
E(T_r)$ is our desired $S$-Steiner tree. Therefore,
$d_G(S)\leq e(T)=s+2t+\sum_{i=3}^rn_i+(n_1-1)+(n_2-1)=n-3$, as
desired.

Next, we consider the case that $x,y$ belong to one connected
component $C_{j_1}$ where $j_1\in \{1,2,\cdots,r\}$. If
$G[V(C_{j_1})\cup \{u\}]$ is $3$-connected, then $G[V(C_{j_1})\cup
\{u\}]\setminus \{x,y\}$ is connected. Therefore, $G[V(C_{j_1})\cup
\{u\}]\setminus \{x,y\}$ contains a spanning tree $T_j$. For each $i
(i\neq j_1, 1\leq i\leq r)$, $G[V(C_{i})\cup \{u\}]$ contains a
spanning tree, say $T_i$. Then the tree $T$ induced by the edges in
$\{uw_i\,|\,1\leq i\leq s\}\cup \{uu_{i}\,|\,1\leq i\leq t\}\cup
\{uv_{i}\,|\,1\leq i\leq t\}\cup E(T_1)\cup E(T_2)\cup \cdots \cup
E(T_r)$ is our desired $S$-Steiner tree. Therefore,
$d_G(S)\leq e(T)=n-3$, as desired.

Suppose $\kappa(G[V(C_{j_1})\cup \{u\}])=2$ and there exists a
vertex $v\in V(C_{j_1})$ such that $\{u,v\}$ is a vertex cut set of
$G[V(C_{j_1})\cup \{u\}]$ and for each component $C_{j_1}^i \
(1\leq i\leq p)$ of $G[V(C_{j_1})\cup \{u\}]\setminus \{u,v\}$, one of the following conditions holds.

\begin{itemize}
\item $G[V(C_{j_1}^i)\cup \{u\}]$ is $3$-connected;

\item $G[V(C_{j_1}^i)\cup
\{v\}]$ is $3$-connected;

\item $\kappa(G[V(C_{j_1}^i)\cup
\{u\}])=\kappa(G[V(C_{j_1}^i)\cup \{v\}])=2$ and $\{z',z''\}$ is not a
common vertex cut set of $G[V(C_{j_1}^i)\cup \{u\}]$ and
$G[V(C_{j_1}^i)\cup \{v\}]$ where $z',z''\in V(C_{j_1}^i)$;

\item
$\kappa(G[V(C_{j_1}^i)\cup \{u\}])=2$ and $\kappa(G[V(C_{j_1}^i)\cup
\{v\}])=1$ and if $\{z',z''\}$ is a vertex cut set of
$G[V(C_{j_1}^i)\cup \{u\}]$ then neither $z'$ nor $z''$ is a cut vertex
of $G[V(C_{j_1}^i)\cup \{v\}]$.
\end{itemize}

It is clear that
$|E_G[C_{j_1}^i),u]|\geq 1$ and $|E_G[C_{j_1}^i),v]|\geq 1$ for each
$i \ (1\leq i\leq p)$.

\begin{itemize}
\item[] Suppose $v\in \{x,y\}$. Without loss of generality, let $v=x$.
Then $y\in \bigcup_{i=1}^pV(C_{j_1}^i)$. Then there exists some
component, say $C_{j_1}^{i_1}$, such that $y\in C_{j_1}^{i_1}$.
Since $G[V(C_{j_1}^{i_1})\cup \{u\}]$ is $2$-connected, it follows
that $G[V(C_{j_1}^{i_1})\cup \{u\}]\setminus y$ is connected and
hence $G[V(C_{j_1}^{i_1})\cup \{u\}]\setminus y$ contains a spanning
tree, say $T_{j_1,i_1}$. For each $i \ (1\leq i\leq p, \ i\neq
i_1)$, because $G[V(C_{j_1}^{i})\cup \{u\}]$ is connected and hence
$G[V(C_{j_1}^{i})\cup \{u\}]$ contains a spanning tree, say
$T_{j_1,i}$. Since $|E_G[C_{j_1}^i),v]|\geq 1$ for each $i \ (1\leq
i\leq p)$, it follows that there exists a component $C_{j_1}^{i_2}$
such that $y\in V(C_{j_1}^{i_2})$ and there exists an edge $vv_p\in
E_G[C_{j_1}^{i_2}),v]$, where $v_p\in V(C_{j_1}^{i_2})$. Set
$T_{j_1}=(\bigcup_{i=1}^pT_{j_1,i})\cup \{vv_p\}$. One can see that
$T_{j_1}$ is a spanning tree of $G[V(C_{j_1})\cup \{u\}]$. For each
$j \ (1\leq i\leq r, \ j\neq j_1)$, because $G[V(C_{j})\cup \{u\}]$
is connected and hence $G[V(C_{j})\cup \{u\}]$ contains a spanning
tree, say $T_{j}$. Then the tree $T$ induced by the edges in
$\{uw_i\,|\,1\leq i\leq s\}\cup \{uu_{i}\,|\,1\leq i\leq t\}\cup
\{uv_{i}\,|\,1\leq i\leq t\}\cup E(T_1)\cup E(T_2)\cup \cdots \cup
E(T_r)$ is our desired $S$-Steiner tree. Therefore,
$d_G(S)\leq e(T)=n-3$, as desired.

\item[] Suppose $v\notin \{x,y\}$. Then $x,y\in \bigcup_{i=1}^pV(C_{j_1}^i)$.
Consider the case that $x,y$ belong to different components in
$\{C_{j_1}^1,C_{j_1}^2,\cdots,C_{j_1}^p\}$. Without loss of
generality, let $x\in V(C_{j_1}^1)$ and $y\in V(C_{j_1}^2)$. Since
$G[V(C_{j_1}^{i})\cup \{u\}] \ (i=1,2)$ is $2$-connected, it follows
that $G[V(C_{j_1}^{1})\cup \{u\}]\setminus x$ contains a spanning
tree $T_{j_1,1}$ and $G[V(C_{j_1}^{2})\cup \{u\}]\setminus y$
contains a spanning tree $T_{j_1,2}$. For each $i \ (3\leq i\leq
p)$, $G[V(C_{j_1}^{i})\cup \{u\}]$ is connected and hence
$G[V(C_{j_1}^{i})\cup \{u\}]$ contains a spanning tree, say
$T_{j_1,i}$. We want to obtain a spanning tree of $G[V(C_{j_1})\cup
\{u\}]\setminus \{x,y\}$ from $T_{j_1,1},T_{j_1,2},\cdots,T_{j_1,p}$ by adding one
edge. If $uv\in E(G)$, then we set
$T_{j_1}=(\bigcup_{i=1}^pT_{j_1,i})\cup \{vu\}$.
Suppose $p\geq 3$.
Since $|E_G[C_{j_1}^i),v]|\geq 1$ for each $i \ (1\leq i\leq p)$, it
follows that there exists a component, say $C_{j_1}^{3}$, such that
$vv_p\in
E_G[C_{j_1}^{3},v]$, where $v_p\in V(C_{j_1}^{3})$. Set
$T_{j_1}=(\bigcup_{i=1}^pT_{j_1,i})\cup \{vv_p\}$.
Suppose $p\geq
2$, and $|E_G[C_{j_1}^1,v]|\geq 2$ or $|E_G[C_{j_1}^2,v]|\geq 2$.
Without loss of generality, let $|E_G[C_{j_1}^1,v]|\geq 2$. Then
there exists two vertices $v_p,v_q\in V(C_{j_1}^{1})$ such that
$vv_p\in E_G[C_{j_1}^{1},v]$ and hence $v_p\neq x$ or $v_q\neq x$.
Without loss of generality, let $v_q\neq x$. Thus, we set
$T_{j_1}=(\bigcup_{i=1}^pT_{j_1,i})\cup \{vv_p\}$. For each $j \
(1\leq i\leq r, \ j\neq j_1)$, because $G[V(C_{j})\cup \{u\}]$ is
connected and hence $G[V(C_{j})\cup \{u\}]$ contains a spanning
tree, say $T_{j}$. Then the tree $T$ induced by the edges in
$\{uw_i\,|\,1\leq i\leq s\}\cup \{uu_{i}\,|\,1\leq i\leq t\}\cup
\{uv_{i}\,|\,1\leq i\leq t\}\cup E(T_1)\cup E(T_2)\cup \cdots \cup
E(T_r)$ is our desired $S$-Steiner tree. Therefore,
$d_G(S)\leq e(T)=n-3$, as desired. Consider the case that $x,y$
belong to same connected component. Without loss of generality, let
$x,y\in V(C_{j_1}^1)$. For $i \ (2\leq i\leq p)$, since
$G[V(C_{j_1}^{i})\cup \{u\}]$ is $2$-connected, it follows that
$G[V(C_{j_1}^{i})\cup \{u\}]$ contains a spanning tree, say
$T_{j_1,i}$. Since $|E_G[C_{j_1}^2,v]|\geq 1$, it follows that
there exists a vertex $v_p\in V(C_{j_1}^{2})$ such that $vv_p\in
E_G[C_{j_1}^{2},v]$. If $G[V(C_{j_1}^{i})\cup \{u\}]$ is
$3$-connected, then $G[V(C_{j_1}^{i})\cup \{u\}]\setminus \{x,y\}$
contains a spanning tree, say $T_{j_1,1}$. If $G[V(C_{j_1}^{i})\cup
\{v\}]$ is $3$-connected, then $G[V(C_{j_1}^{i})\cup \{v\}]\setminus
\{x,y\}$ contains a spanning tree, say $T_{j_1,1}$. Suppose that
$\kappa(G[V(C_{j_1}^1)\cup \{u\}])=\kappa(G[V(C_{j_1}^1)\cup
\{v\}])=2$ and $\{z',z''\}$ is not a common vertex cut set of
$G[V(C_{j_1}^1)\cup \{u\}]$ and $G[V(C_{j_1}^1)\cup \{v\}]$ where
$z',z''\in V(C_{j_1}^1)$. Then $\{x,y\}$ is not a common vertex cut set
of $G[V(C_{j_1}^1)\cup \{u\}]$ and $G[V(C_{j_1}^1)\cup \{v\}]$. If
$\{x,y\}$ is not a vertex cut set of $G[V(C_{j_1}^1)\cup \{u\}]$,
then $G[V(C_{j_1}^1)\cup \{u\}]\setminus \{x,y\}$ contains a
spanning tree $T_{j_1,1}$. If $\{x,y\}$ is not a vertex cut set of
$G[V(C_{j_1}^1)\cup \{v\}]$, then $G[V(C_{j_1}^1)\cup
\{v\}]\setminus \{x,y\}$ contains a spanning tree, say $T_{j_1,1}$.
Suppose that $\kappa(G[V(C_i^j)\cup \{u\}])=2$ and
$\kappa(G[V(C_i^j)\cup \{v\}])=1$ and if $\{z',z''\}$ is a vertex cut
set of $G[V(C_i^j)\cup \{u\}]$ then neither $z'$ nor $z''$ is a cut
vertex of $G[V(C_i^j)\cup \{v\}]$. If $\{x,y\}$ is not a vertex cut
set of $G[V(C_i^j)\cup \{u\}]$, then $G[V(C_{j_1}^1)\cup
\{u\}]\setminus \{x,y\}$ is connected since $G[V(C_{j_1}^1)\cup
\{u\}]$ is $2$-connected. Therefore, $G[V(C_{j_1}^1)\cup
\{u\}]\setminus \{x,y\}$ contains a spanning tree, say $T_{j_1,1}$.
If $\{x,y\}$ is a vertex cut set of $G[V(C_i^j)\cup \{u\}]$, then
neither $x$ nor $y$ is a cut vertex of $G[V(C_i^j)\cup \{v\}]$. Thus
$G[V(C_{j_1}^1)\cup \{v\}]\setminus \{x,y\}$ is connected and hence
$G[V(C_{j_1}^1)\cup \{v\}]\setminus \{x,y\}$ contains a spanning
tree, say $T_{j_1,1}$. Set $T_{j_1}=(\bigcup_{i=1}^pT_{j_1,i})\cup
\{vv_p\}$. For each $j \ (1\leq i\leq r, \ j\neq j_1)$, because
$G[V(C_{j})\cup \{u\}]$ is connected and hence $G[V(C_{j})\cup
\{u\}]$ contains a spanning tree, say $T_{j}$. Then the tree $T$
induced by the edges in $\{uw_i\,|\,1\leq i\leq s\}\cup
\{uu_{i}\,|\,1\leq i\leq t\}\cup \{uv_{i}\,|\,1\leq i\leq t\}\cup
E(T_1)\cup E(T_2)\cup \cdots \cup E(T_r)$ is our desired $S$-Steiner
tree. Therefore, $d_G(S)\leq e(T)=n-3$, as desired.\qed
\end{itemize}
\end{pf}

\begin{lem}\label{lem4}
Let $G$ be a connected graph of order $n$. If $\kappa(G)=2$, then $sdiam_{n-3}(G)=n-3$
if and only if $G$ contains a vertex cut set $\{u,v\}$ and for each connected component $C_i$ of order at least $3$ in
$G\setminus \{u,v\}$, $C_i$ satisfies one of the following two conditions.

$(1)$ $G[V(C_i)\cup \{u,v\}]$ is $3$-connected;

$(2)$ $\kappa(G[V(C_i)\cup \{u,v\}])=2$, both $G[V(C_i)\cup \{u\}]$ and
$G[V(C_i)\cup \{v\}]$ are $2$-connected, for any vertex cut $\{x,y\}\neq \{u,v\}$ of $G[V(C_i)\cup \{u,v\}]$ and any connected component $C_i^j \ (1\leq j\leq s)$ of $G[V(C_i^j)\cup \{u,v\}]\setminus \{x,y\}$, one of the following conditions is true.

$(2.1)$ $G[V(C_i^j)\cup \{x\}]$ is $2$-connected;

$(2.2)$ $G[V(C_i^j)\cup \{y\}]$ is $2$-connected;

$(2.3)$ $\kappa(G[V(C_i^j)\cup \{x\}])=1$, and for any cut vertex $z$ and any component $C_i^{j,k} \ (1\leq k\leq t)$, $|E_G[C_i^{j,k},x]|\geq 1$.

$(2.4)$ $\kappa(G[V(C_i^j)\cup \{y\}])=1$, and for any cut vertex $z$ and any component $C_i^{j,k} \ (1\leq k\leq t)$, $|E_G[C_i^{j,k},x]|\geq 1$.
\end{lem}
\begin{pf}
In one direction, we suppose $sdiam_{n-3}(G)=n-3$. For a vertex cut set $\{u,v\}$, there
exists a connected component $C_{j}$ of order at least $3$ of the
graph $G\setminus \{u,v\}$ such that

\begin{itemize}
\item $\kappa(G[V(C_j)\cup
\{u,v\}])=1$;

\item $\kappa(G[V(C_j)\cup \{u,v\}])=2$ and
$\kappa(G[V(C_i)\cup \{u\}])=1$;

\item $\kappa(G[V(C_j)\cup
\{u,v\}])=2$ and $\kappa(G[V(C_i)\cup \{v\}])=1$;

\item $\kappa(G[V(C_j)\cup
\{u,v\}])=\kappa(G[V(C_j)\cup \{u\}])=\kappa(G[V(C_j)\cup \{v\}])=2$, and there exists a vertex cut set $\{x,y\}$ and a connected component $C_j^{i_1}$ of $G[V(C_j)\cup
\{u,v\}\setminus \{x,y\}$ such that $\kappa(G[V(C_j^{i_1}])\cup \{x\}])=1$, and there exists a cut vertex $z$ and a connected component $C_j^{i_1,k_1}$ such that $|E_G[C_j^{i_1,k_1},x]|=0$.
\end{itemize}

Our aim is to get
a contradiction. Let $H_j=G[V(C_j)\cup
\{u,v\}]$. Suppose that $\kappa(H_j)=2$ and $\kappa(H_j\setminus u)=1$,
or $\kappa(H_j)=2$ and $\kappa(H_j\setminus v)=1$. Without loss of
generality, let $\kappa(H_j)=2$ and $\kappa(H_j\setminus v)=1$.
Since $\kappa(H_j\setminus v)=1$, it follows that there exists a cut
vertex of $H_j\setminus v$, say $x$. One can see that $x\neq u$
since $H_j\setminus \{u,v\}$ is connected. Thus, $H_j\setminus
\{x,v\}$ is disconnected. Let $C_j^1,C_j^2,\cdots,C_j^s \ (s\geq 2)$
be the connected components of $H_j\setminus \{x,v\}$. Clearly,
$u\in \bigcup_{i=1}^sV(C_j^i)$. Without loss of generality, let
$u\in V(C_j^1)$. Then there exists a connected component $C_j^{i_1}
\ (i_1\in \{2,3,\cdots,s\})$ such that $|E_G[u,C_j^{i_1}]|=0$.
Choose $S=V(G)\setminus \{u,v,x\}$. Then $|S|=n-3$. Obviously, any
Steiner tree connecting $S$, say $T$, must occupy the vertices $u$
and $x$, which implies $|V(T)|\geq n-1$ and hence $d_G(S)\geq
e(T)\geq n-2$ and hence $sdiam_{n-3}(G)\geq d_G(S)\geq n-2$, also a
contradiction.

Suppose $\kappa(H_j)=1$. If $v$ is a cut vertex of
$H_j$, then $H_j\setminus v$ is disconnected. Then $u$ belongs to a
connected component $C_j^i$ of $H_j\setminus v$. If $|V(C_j^i)|\geq
2$, then $H_j\setminus \{u,v\}$ is also disconnected, a
contradiction. So we may assume that $V(C_j^i)=\{u\}$. Then $v$ is a
cut vertex of $G$, which contradict to $\kappa(G)=2$. We now suppose
that neither neither $u$ nor $v$ is a cut vertex of $H_j$. Then
there exists a cut vertex of $H_j$, say $x$, such that $x\neq u$ and
$x\neq v$. If $u,v$ belongs to the same connected component of
$H_j\setminus x$, then we choose $S=V(G)\setminus \{u,v,x\}$.
Obviously, $|S|=n-3$ and any Steiner tree connecting $S$, say $T$,
must contain the vertex $x$ and one of $\{u,v\}$, which implies
$|V(T)|\geq n-1$ and hence $d_G(S)\geq e(T)\geq n-2$ and hence
$sdiam_{n-3}(G)\geq d_G(S)\geq n-2$, a contradiction. Suppose $u,v$
belongs to the different connected component of $H_j\setminus x$. One
can also check that $sdiam_{n-3}(G)\geq n-2$ and get a
contradiction.

Suppose that $\kappa(H_j)=\kappa(H_j\setminus \{u\})=\kappa(H_j\setminus \{v\})=2$, and there exists a vertex cut set $\{x,y\}$ and a connected component $C_j^{i_1}$ of $H_j\setminus \{x,y\}$ such that $\kappa(G[V(C_j^{i_1})\cup \{x\}])=1$, and there exists a cut vertex $z$ and a connected component $C_j^{i_1,k_1}$ such that $|E_G[C_j^{i_1,k_1},x]|=0$. Set $S=V(G)\setminus \{x,y,z\}$.
Obviously, $|S|=n-3$ and any Steiner tree connecting $S$, say $T$,
must contain the vertex $z$ and one element $\{x,y\}$, which implies
$|V(T)|\geq n-1$ and hence $d_G(S)\geq e(T)\geq n-2$ and hence
$sdiam_{n-3}(G)\geq d_G(S)\geq n-2$, a contradiction. Suppose $u,v$
belongs to the different connected component of $H_j\setminus x$. One
can also check that $sdiam_{n-3}(G)\geq n-2$ and get a
contradiction.

In another direction, it suffices to show that $d_{G}(S)\leq n-3$ for any
$S\subseteq V(G)$ and $|S|=n-3$. Let $\bar{S}=V(G)\setminus S$.
Clearly, $0\leq |\{u,v\}\cap S|\leq 2$, and each connected component
of $G\setminus \{u,v\}$ is a connected subgraph of order at least
$3$, or an edge of $G$, or an isolated vertex. Let
$w_1,w_2,\cdots,w_s$ be the isolated vertices, $e_1,e_2,\cdots,e_t$
be the independent edges, and $C_1,C_2,\cdots,C_r$ be the connected
components of order at least $3$ in $G\setminus \{u,v\}$. Set
$e_i=u_iv_i \ (1\leq i\leq t)$, $W=\{w_1,w_2,\cdots,w_s\}$,
$U=\{u_1,u_2,\cdots,u_t\}$, $V=\{v_1,v_2,\cdots,v_t\}$. Obviously,
there is at least one edge between $u$ and $\{u_i,v_i\}$, and there
is at least one edge between $v$ and $\{u_i,v_i\}$, and
$uw_i,vw_i\in E(G)$, and $s+2t+\sum_{i=1}^rn_i$ where
$n_i=|V(C_i)|$. Clearly, $0\leq |\bar{S}\cap \{u,v\}|\leq 2$.

Consider the case $|\bar{S}\cap \{u,v\}|=2$. Since $|S|=n-3$, there
exists a vertex $x$ such that $x\notin \{u,v\}$ and $x\notin S$.
Then $x\in W$ or $x\in U\cup V$ or $x\in \bigcup_{i=1}^rV(C_i)$.
\begin{itemize}
\item[] Suppose $x\in W=\{w_1,w_2,\cdots,w_s\}$.
Without loss of generality, let $x=w_1$. Clearly, $G\setminus
\{x,v\}$ is connected. Then $G\setminus \{x,v\}$ contains a spanning
tree $T$, which is a Steiner tree connecting $S$. Observe that
$|V(T)|=n-2$. Therefore, $d_G(S)\leq e(T)\leq n-3$, as desired.

\item[] Suppose $x\in U\cup V$. Without loss of generality, let
$x=v_1$ and $uu_1\in E(G)$. Thus, $G\setminus \{x,v\}$ contains a
spanning tree $T$, which is a Steiner tree connecting $S$. Observe
that $|V(T)|=n-2$. Therefore, $d_G(S)\leq e(T)\leq n-3$, as desired.

\item[] Suppose $x\in \bigcup_{i=1}^rV(C_i)$.
Without loss of generality, let $x\in V(C_1)$. Note that
$G[V(C_1)\cup \{u,v\}]$ is $3$-connected, or $\kappa(G[V(C_1)\cup
\{u,v\}])=2$ and $G[V(C_1)\cup \{u\}]$ is $2$-connected. Then
$G[V(C_1)\cup \{u\}]\setminus x$ is connected and hence
$G[V(C_1)\cup \{u\}]\setminus x$ contains a spanning tree $T_1$. For
each $C_i \ (2\leq i\leq r)$, $G[V(C_i)\cup \{u\}]$ contains a
spanning tree $T_i$. Therefore, the tree $T$ induced by the edges in
$\{uw_1,uw_2,\cdots,uw_s\}\cup \{uu_{1},uu_{2},\cdots, uu_{r}\}\cup
\{u_{1}v_{1},u_{2}v_{2},\cdots,u_{r}v_{r}\}\cup E(T_1)\cup
E(T_2)\cup \cdots \cup E(T_r)$ form a Steiner tree connecting $S$
and so $d_G(S)\leq e(T)=\sum_{i=2}^rn_i+(n_1-1)+2t+s=n-3$, as
desired.
\end{itemize}

Consider the case $|\bar{S}\cap \{u,v\}|=1$. Without loss of
generality, let $v\in \bar{S}$ and $u\notin \bar{S}$. Then $v\notin
S$ and $u\in S$. Since $|S|=n-3$, there exists two vertices $x,y$
such that $x,y\notin \{u,v\}$ and $x,y\notin S$. Then $x,y\in
V(G)\setminus \{u,v\}=W\cup U\cup V \cup (\bigcup_{i=1}^rV(C_i))$.
\begin{itemize}
\item[] Suppose that at least one of $\{x,y\}$ belong to $W$.
Without loss of generality, let $x\in W$ and $x=w_1$. Note that for
each $C_i \ (1\leq i\leq r)$, $G[V(C_i)\cup \{u,v\}]$ is
$3$-connected, or $\kappa(G[V(C_i)\cup \{u,v\}])=2$ and
$G[V(C_i)\cup \{u\}]$ is $2$-connected. Furthermore, $G[V(C_i)\cup
\{u\}]$ contains a spanning tree, say $T_i$. Then the tree $T$
induced by the edges in $\{uw_i\,|\,1\leq i\leq s\}\cup
\{uu_{i}\,|\,1\leq i\leq t\}\cup \{u_{i}v_{i}\,|\,1\leq i\leq
t\}\cup E(T_1)\cup E(T_2)\cup \cdots \cup E(T_r)$ is a Steiner tree
connecting $S$ and so $d_G(S)\leq
e(T)=(s-1)+2t+\sum_{i=1}^rn_i=n-3$, as desired.

\item[] Suppose that at least one of $\{x,y\}$ belong to $\cup_{i=1}^rV(C_i)$.
Without loss of generality, let $x\in \cup_{i=1}^rV(C_i)$. Then
there exists some $C_j$ such that $x\in V(C_j)$. Observe that
$G[V(C_j)\cup \{u\}]\setminus x$ contains a spanning tree, say
$T_j$. For each $C_i \ (i\neq j, 1\leq i\leq r)$, $G[V(C_i)\cup
\{u\}]$ contains a spanning tree, say $T_i$. Then the tree $T$
induced by the edges in $\{uw_i\,|\,1\leq i\leq s\}\cup
\{uu_{i}\,|\,1\leq i\leq t\}\cup \{u_{i}v_{i}\,|\,1\leq i\leq
t\}\cup E(T_1)\cup E(T_2)\cup \cdots \cup E(T_r)$ is a Steiner tree
connecting $S$ and so $d_G(S)\leq e(T)=s+2t+\sum_{i=1,i\neq
j}^rn_i+(n_j-1)=n-3$, as desired.

\item[] Suppose that neither $x$ nor $y$ belong to $W\cup (\cup_{i=1}^rV(C_i))$.
Then $x\in U\cup V$, and $G\setminus \{x,y\}$ is connected and hence
$G\setminus \{x,y\}$ contains a spanning tree $T$. So $d_G(S)\leq
e(T)=n-3$, as desired.
\end{itemize}

Consider the case $|\bar{S}\cap \{u,v\}|=0$. Then $u,v\in S$. Since
$|S|=n-3$, there exists three vertices $x,y,z$ such that
$x,y,z\notin \{u,v\}$ and $x,y,z\notin S$. Then we have $x,y,z\in
V(G)\setminus \{u,v\}=W\cup U\cup V(\cup_{i=1}^rV(C_i))$.
\begin{itemize}
\item[] Suppose that at least two of $\{x,y,z\}$ belong to $W\cup U\cup V$.
Without loss of generality, let $x,y\in W$. Then $G\setminus
\{x,y\}$ is connected and hence $G\setminus \{x,y\}$ contains a
spanning tree $T$. So $d_G(S)\leq e(T)=n-3$, as desired.

\item[] Suppose that only one of $\{x,y,z\}$ belong to $W\cup U\cup V$.
Without loss of generality, let $x\in W\cup U\cup V$. Then $y,z\in
\cup_{i=1}^rV(C_i)$. Then there exists some $C_j$ such that $y\in
V(C_j)$. Since $G[V(C_j)\cup \{u\}]\setminus x$ is $2$-connected, it
follows that $G[V(C_j)\cup \{u\}]\setminus y$ contains a spanning
tree, say $T_j$. Furthermore, for each $C_i \ (i\neq j, 1\leq i\leq
r)$, $G[V(C_i)\cup \{u\}]$ contains a spanning tree, say $T_i$.
Suppose $x\in W$. Without loss of generality, let $x=w_1$. Then the
tree $T$ induced by the edges in $\{uw_i\,|\,2\leq i\leq s\}\cup
\{uu_{i}\,|\,1\leq i\leq t\}\cup \{u_{i}v_{i}\,|\,1\leq i\leq
t\}\cup E(T_1)\cup E(T_2)\cup \cdots \cup E(T_r)$ is our desired
Steiner tree connecting $S$. Suppose $x\in U\cup V$. Without loss of
generality, let $x=u_1$. Then the tree $T$ induced by the edges in
$\{uw_i\,|\,2\leq i\leq s\}\cup \{uu_{i}\,|\,1\leq i\leq t\}\cup
\{u_{i}v_{i}\,|\,2\leq i\leq t\}\cup E(T_1)\cup E(T_2)\cup \cdots
\cup E(T_r)$ is our desired Steiner tree connecting $S$. Therefore,
$d_G(S)\leq e(T)=n-3$, as desired.

\item[] Suppose that none of $\{x,y,z\}$ belong to $U\cup V\cup W$.
Then $x,y,z\in \cup_{i=1}^rV(C_i)$. Consider the case that $x,y,z$
belong to three connected components, say $C_1,C_2,C_3$. Without loss of generality, let $x\in V(C_1)$, $y\in V(C_2)$ and $z\in V(C_3)$. Since $G[V(C_i)\cup \{u\}] \ (i=1,2,3)$ is $2$-connected, it follows that $G[V(C_1)\cup \{u\}]\setminus x$, $G[V(C_2)\cup \{u\}]\setminus y$ and $G[V(C_3)\cup \{u\}]\setminus z$ are all connected. Therefore, $G[V(C_1)\cup \{u\}]\setminus x$ contains a spanning tree, say $T_1$; $G[V(C_2)\cup \{u\}]\setminus y$ contains a spanning tree, say $T_2$; $G[V(C_3)\cup \{u\}]\setminus z$ contains a spanning tree, say $T_3$. Furthermore, for each $C_i \ (4\leq i\leq r)$,
$G[V(C_i)\cup \{u\}]$ contains a spanning tree, say $T_i$. Then the
tree $T$ induced by the edges in $E(T_1)\cup E(T_2)\cup \cdots \cup E(T_r)$ is our desired
Steiner tree connecting $S$. Therefore, $d_G(S)\leq
e(T)=\sum_{i=1}^rn_i-3=n-5$, as desired. Consider the case that $x,y,z$
belong to two connected components, say $C_1,C_2$. Without loss of generality, let $x,y\in V(C_1)$ and $z\in V(C_3)$. Since $G[V(C_i)\cup \{u\}] \ (i=1,2)$ is $2$-connected, it follows that $G[V(C_1)\cup \{u\}]\setminus x$ and $G[V(C_2)\cup \{u\}]\setminus z$ are both connected. Therefore, $G[V(C_1)\cup \{u\}]\setminus x$ contains a spanning tree, say $T_1$, and $G[V(C_2)\cup \{u\}]\setminus z$ contains a spanning tree, say $T_2$. Furthermore, for each $C_i \ (3\leq i\leq r)$,
$G[V(C_i)\cup \{u\}]$ contains a spanning tree, say $T_i$. Then the
tree $T$ induced by the edges in $E(T_1)\cup E(T_2)\cup \cdots \cup E(T_r)$ is our desired $S$-Steiner tree. Therefore, $d_G(S)\leq
e(T)=\sum_{i=1}^rn_i-2=n-4$, as desired. Consider the case that $x,y,z$
belong to the same connected component, say $C_1$. If $\{x,y\}$ or $\{y,z\}$ or $\{x,z\}$ is not a vertex cut set of $G[V(C_1)\cup \{u,v\}]$, then
$G[V(C_1)\cup \{u,v\}]\setminus \{x,y\}$ is connected, it follows that
$G[V(C_1)\cup \{u,v\}]\setminus \{x,y\}$ contains a spanning tree, say
$T_1$. Furthermore, for each $C_i \ (2\leq i\leq r)$,
$G[V(C_i)\cup \{u\}]$ contains a spanning tree, say $T_i$. Then the
tree $T$ induced by the edges in $E(T_1)\cup E(T_2)\cup \cdots \cup E(T_r)$ is our desired
Steiner tree connecting $S$. Therefore, $d_G(S)\leq
e(T)=\sum_{i=1}^rn_i-2+2-1=n-3$, as desired.

Suppose that $\{x,y\}$, $\{y,z\}$ and $\{x,z\}$ are all vertex cut sets of $G[V(C_1)\cup \{u,v\}]$. We consider the connected component $C_1^{j} \ (1\leq j\leq s)$ of $G[V(C_1)\cup \{u,v\}]\setminus \{x,y\}$.
If $G[V(C_1^{j})\cup \{x\}]$ is $2$-connected, then $G[V(C_1^{j})\cup \{x\}]\setminus \{x,v\}=G[V(C_1^{j})\cup \{u\}]\setminus \{x\}$ is connected, then
$G[V(C_1^{j})\cup \{x\}]\setminus \{z\}$ contains a spanning tree, say
$T_{1,j}$. Furthermore, for each $C_1^{i} \ (1\leq i\leq s, i\neq j)$,
$G[V(C_1^{i})\cup \{y\}]$ contains a spanning tree, say $T_{1,i}$. Let $T_{1}$ denote the tree induced by the edges in $E(T_{1,1})\cup E(T_{1,2})\cup \cdots \cup E(T_{1,r})$. Furthermore, for each $C_i \ (2\leq i\leq r)$,
$G[V(C_i)\cup \{u\}]$ contains a spanning tree, say $T_i$. Then the
tree $T$ induced by the edges in $E(T_1)\cup E(T_2)\cup \cdots \cup E(T_r)$ is our desired $S$-Steiner tree. Therefore, $d_G(S)\leq
e(T)=\sum_{i=1}^rn_i-1=n-3$, as desired.

Suppose that $\kappa(G[V(C_i^j)\cup \{x\}])=1$ and for any cut vertex $z$ and any component $C_i^{j,k} \ (1\leq k\leq t)$, $|E_G[C_i^{j,k},x]|\geq 1$. Then $G[V(C_i^j)\cup \{x\}]\setminus \{z\}$ is connected, and hence $G[V(C_i^j)\cup \{x\}]\setminus \{z\}$ contains a spanning tree, say , say
$T_{1,j}$. Furthermore, for each $C_1^{i} \ (1\leq i\leq s, i\neq j)$,
$G[V(C_1^{i})\cup \{x\}]$ contains a spanning tree, say $T_{1,i}$. Let $T_{1}$ denote the tree induced by the edges in $E(T_{1,1})\cup E(T_{1,2})\cup \cdots \cup E(T_{1,r})$. Furthermore, for each $C_i \ (2\leq i\leq r)$,
$G[V(C_i)\cup \{u\}]$ contains a spanning tree, say $T_i$. Then the
tree $T$ induced by the edges in $E(T_1)\cup E(T_2)\cup \cdots \cup E(T_r)$ is our desired $S$-Steiner tree. Therefore, $d_G(S)\leq
e(T)=\sum_{i=1}^rn_i-1=n-3$, as desired.
\end{itemize}

From the above argument, $d_{G}(S)\leq n-3$ for any $S\subseteq
V(G)$ and $|S|=n-3$. So $sdiam_{n-3}(G)\leq n-3$. \qed
\end{pf}

\begin{pro}\label{pro3}
Let $G$ be a connected graph of order $n$. Then $sdiam_{n-3}(G)=n-3$
if and only if $G$ satisfies one of the following conditions.

$(1)$ $\kappa(G)=3$;

$(2)$ $\kappa(G)=2$ and $G$ contains a vertex cut set $\{u,v\}$ and
for each connected component $C_i$ of order at least $3$ in
$G\setminus \{u,v\}$, $C_i$ satisfies one of the following two
conditions.

$(2.1)$ $G[V(C_i)\cup \{u,v\}]$ is $3$-connected;

$(2.2)$ $\kappa(G[V(C_i)\cup \{u,v\}])=2$, both $G[V(C_i)\cup
\{u\}]$ and $G[V(C_i)\cup \{v\}]$ are $2$-connected, for any vertex
cut $\{x,y\}\neq \{u,v\}$ of $G[V(C_i)\cup \{u,v\}]$ and any
connected component $C_i^j \ (1\leq j\leq s)$ of $G[V(C_i^j)\cup
\{u,v\}]\setminus \{x,y\}$, one of the following conditions is true.

$(2.2.1)$ $G[V(C_i^j)\cup \{x\}]$ is $2$-connected;

$(2.2.2)$ $G[V(C_i^j)\cup \{y\}]$ is $2$-connected;

$(2.2.3)$ $\kappa(G[V(C_i^j)\cup \{x\}])=1$, and for any cut vertex
$z$ and any component $C_i^{j,k} \ (1\leq k\leq t)$,
$|E_G[C_i^{j,k},x]|\geq 1$.

$(2.2.4)$ $\kappa(G[V(C_i^j)\cup \{y\}])=1$, and for any cut vertex
$z$ and any component $C_i^{j,k} \ (1\leq k\leq t)$,
$|E_G[C_i^{j,k},x]|\geq 1$.

$(3)$ $G$ contains only one cut vertex $u$; for each connected
component $C_i$ of order at least $3$ in $G\setminus u$,
$G[V(C_i)\cup \{u\}]$ is $3$-connected, or $\kappa(G[V(C_i)\cup
\{u\}])=2$ and there exists a vertex $v\in V(C_i)$ such that
$\{u,v\}$ is a vertex cut set of $G[V(C_i)\cup \{u\}]$, and for each
component $C_i^j \ (1\leq j\leq p)$ of $G[V(C_i)\cup \{u\}]\setminus
\{u,v\}$, one of the following conditions holds:

$\bullet$ $uv\in E(G)$;

$\bullet$ $p\geq 3$;

$\bullet$ $p=2$, and $|E_G[v,V(C_i^1)]|\geq 2$ or
$|E_G[v,V(C_i^2)]|\geq 2$.

and one of the following conditions holds:

$\bullet$ $G[V(C_i^j)\cup \{u\}]$ is $3$-connected;

$\bullet$ $G[V(C_i^j)\cup \{v\}]$ is $3$-connected;

$\bullet$ $\kappa(G[V(C_i^j)\cup \{u\}])=\kappa(G[V(C_i^j)\cup
\{v\}])=2$ and $\{y,z\}$ is not a common vertex cut set of
$G[V(C_i^j)\cup \{u\}]$ and $G[V(C_i^j)\cup \{v\}]$ where $z',z''\in
V(C_i^j)$;

$\bullet$ $\kappa(G[V(C_i^j)\cup \{u\}])=2$ and
$\kappa(G[V(C_i^j)\cup \{v\}])=1$ and if $\{z',z''\}$ is a vertex
cut set of $G[V(C_i^j)\cup \{u\}]$, then neither $z'$ nor $z''$ is a
cut vertex of $G[V(C_i^j)\cup \{v\}]$.
\end{pro}
\begin{pf}
In one direction, we suppose $sdiam_{n-3}(G)=n-3$. From Lemma
\ref{lem2}, we have $\kappa(G)\leq 3$. If $\kappa(G)=2$, then the
result holds by Lemma \ref{lem4}. If $\kappa(G)=1$, then the result
holds by Lemma \ref{lem3}.

Conversely, we suppose that $G$ satisfies the conditions of this
theorem. From Lemma \ref{lem2}, we know that $sdiam_{n-3}(G)\geq
n-3$. So it suffices to show that $sdiam_{n-3}(G)\leq n-3$. If
$\kappa(G)=3$, then $sdiam_{n-3}(G)\leq sdiam_{n-2}(G)=n-3$ by $(2)$
of Observation \ref{obs1} and $(1)$ of Lemma \ref{lem2}. From $(1)$
of this theorem, we have $sdiam_{n-3}(G)\geq n-3$. Therefore,
$sdiam_{n-3}(G)=n-3$. If $\kappa(G)=2$ and $G$ satisfies Condition
$(2)$, then $sdiam_{n-3}(G)\leq n-3$ by Lemma \ref{lem4}. If
$\kappa(G)=1$ and $G$ satisfies Condition $(3)$, then
$sdiam_{n-3}(G)\leq n-3$ by Lemma \ref{lem3}. \qed
\end{pf}






\begin{cor}\label{cor3}
Let $G$ be a connected graph of order $n$. Then $sdiam_{n-3}(G)=n-2$
if and only if $G$ satisfies one of the following conditions.

$(1)$ $\kappa(G)=2$; for a vertex cut set $\{u,v\}$, there exists a
connected component $C_{j}$ of order at least $3$ of the graph
$G\setminus \{u,v\}$ such that

\begin{itemize}
\item $\kappa(G[V(C_j)\cup
\{u,v\}])=1$;

\item $\kappa(G[V(C_j)\cup \{u,v\}])=2$ and
$\kappa(G[V(C_i)\cup \{u\}])=1$;

\item $\kappa(G[V(C_j)\cup
\{u,v\}])=2$ and $\kappa(G[V(C_i)\cup \{v\}])=1$;

\item $\kappa(G[V(C_j)\cup
\{u,v\}])=\kappa(G[V(C_j)\cup \{u\}])=\kappa(G[V(C_j)\cup
\{v\}])=2$, and there exists a vertex cut set $\{x,y\}$ and a
connected component $C_j^{i_1}$ of $G[V(C_j)\cup \{u,v\}\setminus
\{x,y\}$ such that $\kappa(G[V(C_j^{i_1}])\cup \{x\}])=1$, and there
exists a cut vertex $z$ and a connected component $C_j^{i_1,k_1}$
such that $|E_G[C_j^{i_1,k_1},x]|=0$.
\end{itemize}

$(2)$ there exist exactly two cut vertices in $G$.

$(3)$ $G$ contains only one cut vertex $u$ such that there
exists a connected component $C_j$ of order at least $3$ in
$G\setminus u$ satisfying one of the following.

$(3.1)$ $\kappa(G[V(C_j)\cup \{u\}])=1$;

$(3.2)$ $\kappa(G[V(C_j)\cup \{u\}])=2$ and $\{v,u\}$
is not a vertex cut set of $G[V(C_j)\cup \{u\}]$ for any $v\in V(C_j)$;

$(3.3)$ $\kappa(G[V(C_j)\cup \{u\}])=2$, $\{v,u\}$ is a
vertex cut set of $G[V(C_j)\cup \{u\}]$, and there exists a component
$C_j^{i'}$ of $G[V(C_j)\cup \{u\}]\setminus \{u,v\}$ satisfying one
of the following conditions.

$\bullet$ $\kappa(G[V(C_j^{i'}\cup \{u\}])=1$.

$\bullet$ $\kappa(G[V(C_j^{i'}\cup \{u\}])=\kappa(G[V(C_j^{i'}\cup
\{v\}])=2$ and $\{y,z\}$ is a common vertex cut set of
$G[V(C_j^{i'}\cup \{u\}]$ and $V(C_j^{i'}\cup \{v\}]$ for any
$y,z\in V(C_i^j)$.

$\bullet$ $\kappa(G[V(C_j^{i'}\cup \{u\}])=2$ and
$\kappa(G[V(C_j^{i'}\cup \{v\}])=1$ and if $\{y,z\}$ is a common
vertex cut set of $G[V(C_j^{i'}\cup \{u\}]$ where $y,z\in V(C_i^j)$,
then at least one of $\{y,z\}$ is a cut vertex of $V(C_j^{i'}\cup
\{v\}]$.

$\bullet$ $uv\notin E(G)$, $p=2$ and there is only one edge between
$v$ and each connected component of $G[V(C_j)\cup \{u\}]\setminus
\{u,v\}$.
\end{cor}

From the above lemmas and corollary, we conclude the following
theorem.

\begin{thm}\label{th3}
Let $G$ be a connected graph of order $n$. Then

$(1)$ $sdiam_{n-3}(G)=n-4$ if and only if $\kappa(G)\geq 4$.

$(2)$ $sdiam_{n-3}(G)=n-3$ if and only if $G$ satisfies the
conditions of Proposition \ref{pro3}.

$(3)$ $sdiam_{n-3}(G)=n-2$ if and only if $G$ satisfies the
conditions of Corollary \ref{cor3}.

$(4)$ $sdiam_{n-3}(G)=n-1$ if and only if there are at least three
cut vertices in $G$.
\end{thm}

\end{document}